\documentclass[12pt]{article}
\usepackage[T2A]{fontenc}
\usepackage[english]{babel}
\usepackage[tbtags]{amsmath}
\usepackage{amsfonts,amssymb}
\usepackage{amssymb}
\usepackage{graphicx}
\usepackage{amscd}

\usepackage{fullpage}
\usepackage{float}

\usepackage{graphics,amsmath,amssymb}
\usepackage{amsthm}
\usepackage{amsfonts}
\usepackage{mathrsfs}

\title{Application of a composition of generating functions for obtaining explicit formulas of polynomials}

\author{D.~V.~Kruchinin, V.~V.~Kruchinin\\
\small Tomsk State University of Control Systems and Radioelectronics, Russian Federation\\
\small \texttt{kruchininDm@gmail.com}\\
}

\begin{document}
\maketitle

\begin{abstract}
Using notions of composita and composition of generating functions we obtain explicit formulas for Chebyshev polynomials, Legendre polynomials, Gegenbauer polynomials, Associated Laguerre polynomials,  Stirling polynomials, Abel polynomials, Bernoulli Polynomials of the Second Kind, Generalized Bernoulli polynomials, Euler Polynomials, Peters polynomials, Narumi polynomials, Humbert polynomials, Lerch polynomials and Mahler polynomials. 

Key words: Composita; Generating Function; Composition of Generating Function; Chebyshev polynomials; Legendre polynomials; Gegenbauer polynomials; Associated Laguerre polynomials;  Stirling polynomials; Abel polynomials; Bernoulli Polynomials of the Second Kind; Generalized Bernoulli polynomials; Euler Polynomials; Peters polynomials; Narumi polynomials; Humbert polynomials; Lerch polynomials; Mahler polynomials.
\end{abstract}

\theoremstyle{plain}
\newtheorem{theorem}{Theorem}
\newtheorem{corollary}[theorem]{Corollary}
\newtheorem{lemma}[theorem]{Lemma}
\newtheorem{proposition}[theorem]{Proposition}

\theoremstyle{definition}
\newtheorem{definition}[theorem]{Definition}
\newtheorem{example}[theorem]{Example}
\newtheorem{conjecture}[theorem]{Conjecture}
\theoremstyle{remark}
\newtheorem{remark}[theorem]{Remark}

\newtheorem{Theorem}{Theorem}[section]
\newtheorem{Proposition}[Theorem]{Proposition}
\newtheorem{Corollary}[Theorem]{Corollary}

\theoremstyle{definition}
\newtheorem{Example}[Theorem]{Example}
\newtheorem{Remark}[Theorem]{Remark}
\newtheorem{Problem}[Theorem]{Problem}
\newtheorem{state}[Theorem]{Statement}
\makeatletter
\def\rdots{\mathinner{\mkern1mu\raise\p@
\vbox{\kern7\p@\hbox{.}}\mkern2mu
\raise4\p@\hbox{.}\mkern2mu\raise7\p@\hbox{.}\mkern1mu}}
\makeatother

\section{Introduction}
A polynomial is a mathematical expression involving a sum of powers in one or more variables multiplied by coefficients.
The use of polynomials in other areas of mathematics is very impressive: continued fractions, operator theory, analytic functions, interpolation, approximation theory, numerical analysis, electrostatics, statistical quantum mechanics, special functions, number theory, combinatorics, stochastic processes, sorting and data compression, and etc.

In this paper we use the method of obtaining expressions for polynomials based on the composition of generating functions, which was presented by the authors in \cite{ICNAAM2012}.

The generating functions have important role in many branches of Mathematics
and Mathematical Physics. In the literature, one can find extensive investigations related to the generating functions for many polynomials \cite{Roman1984,Simsek}. 

In this paper we obtain explicit formulas for Chebyshev polynomials, Legendre polynomials, Gegenbauer polynomials, Associated Laguerre polynomials,  Stirling polynomials, Abel polynomials, Bernoulli Polynomials of the Second Kind, Generalized Bernoulli polynomials, Euler Polynomials, Peters polynomials, Narumi polynomials, Humbert polynomials, Lerch polynomials and Mahler polynomials.
\section{Preliminary}

In the paper \cite{KruCompositae} authors introduced the notion of the \textit{composita} of a given ordinary generating function $F(t)=\sum_{n>0}g(n)t^n$.

Suppose $F(t) = \sum_{n>0} f(n) t^n$ is the generating function, in which there is no free term $f(0)=0$.
From this generating function we can write the following condition:
\begin{equation}
[F(t)]^k=\sum_{n>0} F(n,k)t^n.
\end{equation}
The expression $F(n,k)$ is the \textit{composita} and it's denoted by $F^{\Delta}(n, k) $. 

Also the \textit{composita} can be written based on the composition of $n$:\\
The \textit{composita} is the function of two variables defined by
\begin{equation}
\label{Fnk0}F^{\Delta}(n,k)=\sum_{\pi_k \in C_n}{f(\lambda_1)f(\lambda_2)\ldots f(\lambda_k)},
\end{equation}
where $C_n$ is a set of all compositions of an integer $n$, $\pi_k$ is the composition $\sum_{i=1}^k\lambda_i=n$ into $k$ parts exactly.

Louis Comtet\cite[ p.\ 141]{Comtet_1974} considered similar objects and identities for exponential generating functions, and called them potential polynomials. In this paper we consider the more general case of generating functions: ordinary generating functions.

The expression $ F^{\Delta}(n,k)$ takes a triangular form
$$
\begin{array}{ccccccccccc}
&&&&& F_{1,1}^{\Delta}\\
&&&& F_{2,1}^{\Delta} && F_{2,2}^{\Delta}\\
&&& F_{3,1}^{\Delta} && F_{3,2}^{\Delta} && F_{3,3}^{\Delta}\\
&& F_{4,1}^{\Delta} && F_{4,2}^{\Delta} && F_{4,3}^{\Delta} && F_{4,4}^{\Delta}\\
& \vdots && \vdots && \vdots && \vdots && \vdots\\
F_{n,1}^{\Delta} && F_{n,2}^{\Delta} && \ldots && \ldots && F_{n,n-1}^{\Delta} && F_{n,n}^{\Delta}\\
\end{array}
$$

For instance, we obtain the composita of the generating function $G(t,a,b)=at+bt^2$.

Raising this generating function to the power of $k$ and applying the binomial theorem, we obtain 

$$[G(t,a,b)]^k=t^k(a+bt)^k=x^k\sum_{m=0}^k \binom{k}{m}a^{k-m}b^mt^m. 
$$

Substituting $n$ for $m+k$, we get the following expression
$$[G(t,a,b)]^k=\sum_{n=k}^{2k} \binom{k}{n-k}a^{2k-n}b^{n-k}t^n=\sum_{n=k}^{2k}G^{\Delta}(n,k,a,b)t^n. 
$$
Therefore, the composita is 
\begin{equation}
\label{Gnkab}
G^{\Delta}(n,k,a,b)=\binom{k}{n-k}a^{2k-n}b^{n-k}. 
\end{equation}
In the table \ref{tab:a} we show a few compositae of known generating functions \cite{Comtet_1974, Wilf_1994}. 

\begin{table}[h]
\begin{center}
\setlength\arrayrulewidth{1pt}
\renewcommand{\arraystretch}{1,3}
\begin{tabular}{|ccc|ccc|}
\hline
&\textbf{Generating function $G(x)$}& & & \textbf{Composita $G^{\Delta}(n,k)$}&\\
\hline
&$at+bt^2$ &&& $a^{2k-n}b^{n-k}\binom{k}{n-k}$&\\ \hline
&$\frac{bt}{1-ax}$ &&&$\binom{n-1}{k-1}a^{n-k}b^k$&\\  \hline
&$\ln(1+t)$ &&& $\frac{k!}{n!}\genfrac{[}{]}{0pt}{}{n}{k}$& \\ \hline
&$e^t-1$  &&& $\frac{k!}{n!}\genfrac{\{}{\}}{0pt}{}{n}{k}$& \\ \hline
\end{tabular}
\caption{Examples of generating functions and their compositae}
\label{tab:a}
\end{center}
\end{table}

Here $\genfrac{[}{]}{0pt}{}{n}{k}$ and $\genfrac{\{}{\}}{0pt}{}{n}{k}$ stand for the Stirling numbers of the first kind and of the second kind, respectively (see \cite{Comtet_1974,ConcreteMath}).

The Stirling numbers of the first kind $\genfrac{[}{]}{0pt}{}{n}{k}$ count the number of permutations of $n$ elements with $k$ disjoint cycles. The Stirling numbers of the first kind are defined by the following generating function
$$
\psi_k(x)=\sum_{n\geq k}  \genfrac{[}{]}{0pt}{}{n}{k}  \frac{x^n}{n!}=\frac{1}{k!}\ln^k(1+x).
$$

The Stirling numbers of the second kind $\genfrac{\{}{\}}{0pt}{}{n}{k}$ count the number of ways to partition a set of $n$ elements into $k$ nonempty subsets. A general formula for the Stirling numbers of the second kind is given as follows:
$$\genfrac{\{}{\}}{0pt}{}{n}{k}=\frac{1}{k!}\sum_{j=0}^k(-1)^{k-j}\binom{k}{j}j^n.$$

The Stirling numbers of the second kind are defined by the following generating function
$$
\Phi_k(x)=\sum_{n\geq k} \genfrac{\{}{\}}{0pt}{}{n}{k} \frac{x^n}{n!}=\frac{1}{k!}(e^x-1)^k.
$$

Next we show some operations with compositae.

\begin{itemize}
\item Suppose $F(t)=\sum_{n>0} f(n)t^n$ is the generating function, and $F^{\Delta}(n,k)$ is the composita of this generating function, and $\alpha$ is constant.

For the generating function $A(t)=\alpha F(t)$  the composita is equal to 
\begin{equation}
A^{\Delta}(n,k)=\alpha^k F^{\Delta}(n,k).
\end{equation}

\item Suppose $F(t)=\sum_{n>0} f(n)t^n$ is the generating function, and $F^{\Delta}(n,k)$ is the composita of this generating function, and $\alpha$ is constant.

For the generating function $A(t)=F(\alpha t)$ the composita is equal to
\begin{equation}
A^{\Delta}(n,k)=\alpha^n F^{\Delta}(n,k).
\end{equation}

\item Suppose we have the generating function $F(t)=\sum_{n>0} f(n)t^n$ , composita of this generating function $F^{\Delta}(n,k)$; the following generating functions $B(t)=\sum_{n\geqslant 0} b(n)x^n$  and $[B(t)^k]=\sum_{n\geqslant 0}B(n,k)t^n$. 
Then for the generating function $A(t)=F(t)B(t)$  the composita is equal to 
\begin{equation}
A^{\Delta}(n,k)=\sum_{i=k}^{n} F^{\Delta}(i,k)B(n-i,k).
\end{equation}

\item Suppose we have the generating functions $F(t)=\sum_{n>0} f(n)t^n$, $G(t)=\sum_{n>0} g(n)t^n$, and their compositae $F^{\Delta}(n,k)$ , $G^{\Delta}(n,k)$ accordingly.
Then for the generating function $A(t)=F(t)+G(t)$  the composita is equal to
\begin{equation}
A^{\Delta}(n,k)=F^{\Delta}(n,k)+\sum_{j=1}^{k-1}\binom{k}{j}\sum_{i=j}^{n-k+j}F^{\Delta}(i,j)G^{\Delta}(n-i,k-j)+G^{\Delta}(n,k).
\end{equation}

\item Suppose we have the generating functions $R(t)=\sum_{n>0} r(n)t^n$, $F(t)=\sum_{n>0} f(n)t^n$, and their compositae $F^{\Delta}(n,k)$ , $R^{\Delta}(n,k)$ accordingly.
Then for the composition of generating functions $A(t)=F\left(R(t)\right)$  the composita is equal to 
\begin{equation}
\label{compCompositon}
A^{\Delta}(n,m)=\sum_{k=m}^n F^{\Delta}(n,k)R^{\Delta}(k,m).
\end{equation}

\item For the composition of generating functions $ A(t)=R(F(t))$ the following condition holds:
\begin{equation}
\label{composition}
a(n)=\sum_{k=1}^n F^{\Delta}(n,k)r(k),\qquad a(0)=r(0),
\end{equation}
where $A(t)=\sum_{n\geq 0} a(n)t^n.$
\end{itemize}

Let us consider the generating function of two variables $f(x,t)=2xt-t^2$.
Then 
$F(x,t)^k=\sum_{n\geq k} F^{\Delta}(n,k,x)t^n$, where, according to table \ref{tab:a},
$$
F^{\Delta}(n,k,x)={\binom{k}{n-k}}\,2^{2\,k-n}\,\left(-1\right)^{n-k}\,x^{2\,k-n}.
$$
The triangular form of this composita is the following

$$\begin{array}{lllllll}2\,x\\-1&4\,x^2\\0&-4\,x&8\,x^3\\0&1&-12\,x^2&16\,x^4\\0&0&6\,x&-32\,x^3&32\,x^5\\0&0&-1&24\,x^2&-80\,x^4&64\,x^6\\0&0&0&-8\,x&80\,x^3&-192\,x^5&128\,x^7\\\end{array}$$

Now let's consider a generating function $R(t)=\sum_{n \geq 0} r(n)t^n $. Then the composition of generating functions $P(t,x)=R(2xt-t^2)$ defines a family of polynomials, which have the following form:
\begin{equation}
\label{P_n_common}
P_n(x)=\sum_{k=1}^n F^{\Delta}(n,k,x)r(k)=\sum_{k=\left \lceil {{n}\over{2}} \right \rceil}^n{\binom{k}{n-k}}\,2^{2\,k-n}\,\left(-1\right)^{n-k}\,x^{2\,k-n}r(k),\quad P_0(x)=r(0).
\end{equation}

With this we can consider an application of this approach to obtain closed formulas for Chebyshev polynomials, Legendre polynomials, Gegenbauer polynomials, Associated Laguerre polynomials,  Stirling polynomials, Abel polynomials, Bernoulli Polynomials of the Second Kind, Generalized Bernoulli polynomials, Euler Polynomials, Peters polynomials, Narumi polynomials, Humbert polynomials, Lerch polynomials and Mahler polynomials.

\section{Chebyshev Polynomials}
 
We start by considering Chebyshev Polynomials.
Chebyshev polynomials \cite{Chebyshev} are a sequence of orthogonal polynomials which are related to de Moivre's formula and which can be defined recursively. One usually distinguishes between Chebyshev polynomials of the first kind which are denoted $T_n$ and Chebyshev polynomials of the second kind which are denoted $U_n$.

Chebyshev polynomials form a special class of polynomials especially suited for
approximating other functions. They are widely used in many areas of numerical analysis:
uniform approximation, least-squares approximation, numerical solution of ordinary and
partial differential equations (the so-called spectral or pseudospectral methods), and so on \cite{NumMeth}.

The generating function for the Chebyschev polynomials  of the first kind is given as follows \cite{MathWolfram}:
$$
\frac{1-tx}{1-2tx+t^2}=\sum_{n \geq 0} T_n(x)t^n.
$$

We can represent this expression as follows:
$$
\frac{1-tx}{1-2tx+t^2}=\left[\frac{1}{2}\ln\left(\frac{1}{1-f(x,t)}\right)\right]'.
$$

Therefore, according to (\ref{P_n_common}), the expression of the polynomial $T_n(x)$ has following form:

$$
T_n(x)=n\sum_{k=\left \lceil {{n}\over{2}} \right \rceil}^n {k \choose n-k}\frac{2^{2k-n-1}}{k}\left(-1\right)^{n-k}\,x^{2\,k-n}.
$$ 

The Chebyshev polynomials of the second kind are defined by the  generating function \cite{MathWolfram}:

$$
\frac{1}{1-2tx+t^2}=\sum_{n\geq 0} U_n(x) t^n.
$$

Therefore, the expression of the composition generating function $\frac{1}{1-f(x,t)}$ defines $U_n(x)$
$$
U_n(x)=\sum_{k=\left \lceil {{n}\over{2}} \right \rceil}^{n}{{{k}\choose{n-k}}\,2^{2\,k-n}\,\left(-1\right)^{n-k}
 \,x^{2\,k-n}}.
$$

\section{Legendre Polynomials}

Legendre polynomials $P_n(x)$ are solutions to Legendre's differential equation:
$$
\frac{d}{dx}\left((1-x^2)\frac{d}{dx}P_n(x)\right)+n(n+1)P_n(x)=0.
$$

Legendre polynomials are important in problems involving spheres or spherical coordinates. Due to their
orthogonality properties they are also useful in numerical analysis.

The generating function for the Legendre polynomials is given as follows \cite{MathWolfram}:
$$
\frac{1}{\sqrt{1-2xt+t^2}}=\sum_{n\geq 0} P_n(x)t^n.
$$

In physics, this generating function is the basis for multipole expansions \cite{ArfWeb}.

For obtaining $P_n(x)$ we consider the generating function as a composition of generating functions
$$
\frac{1}{\sqrt{1-f(x,t)}}.
$$

For the generating function $ R(x)=\frac{1}{\sqrt{1-t}} $ the coefficients are determined by the expression:
$r(n)={\frac{1}{4^n}{2n\choose n}}$.

Then, according to (\ref{P_n_common}), we obtain
$$P_n(x)=\frac{1}{2^n}{{\sum_{k=\left \lceil {{n}\over{2}} \right \rceil}^{n}{{{k
 }\choose{n-k}}\,{{{2\,k}\choose{k}}\,\left(-1\right)^{n-k}\,x^{2\,k-n
 }}}}}.$$

\section{Gegenbauer Polynomials}

The Gegenbauer polynomials $ C_n^{(\alpha)} $ are solutions to the Gegenbauer differential equation for integer $n$.
$$
(1-x^2)y''-(2\alpha+1)xy'+n(n+2\alpha)y=0
$$
They generalize Legendre polynomials and Chebyshev polynomials.

Legendre polynomials are important in applied mathematics, in the implementation of spectral and pseudo-spectral methods.
Gegenbauer polynomials are a convenient basis for polynomial approximations since they are eigenfunctions
of corresponding differential operators.

The generating function for the Gegenbauer polynomials  is given as follows \cite{SteWei}:

$$
\left({1-2xt+t^2}\right)^{-\alpha}=\sum C_n^{(\alpha)} t^n.
$$ 
Let us consider the generating function as a composition of generating functions
$$
\left({1-2xt+t^2}\right)^{-\alpha}=e^{{\alpha}\ln\left(\frac{1}{1-f(x,t)}\right)}
$$ 
According to (\ref{compCompositon}), the composita of the generating function $h(x)=\ln\left(\frac{1}{1-f(x,t)}\right)$ is given as follows
$$H^{\Delta}(n,m)=m!\,\left(-1\right)^{n-m}\,\sum_{k=m}^{n}{{{\left[{k \atop m}\right]\,{{k}\choose{n-k}}\,2^{2\,k-n}\,x^{2\,k-n}}\over{k!}}}$$
Therefore, using (\ref{composition}), we obtain
$$C_n^{(\alpha)}(x)=\sum_{m=0}^{n}{\alpha^{m}\,\left(-1\right)^{n-m}\,\sum_{k={\it max} \left(m , \left \lceil {{n}\over{2}} \right \rceil\right)}^{n}{{{\frac{2^{2k-n}}{k!} \left[{k \atop m}\right]\,{{k}\choose{n-k}}\,\,x
 ^{2\,k-n}}}}}.$$
 
Let us consider other composition of generating functions
$$
\left[\frac{1}{1-f(t,x)}\right]^{\alpha}.
$$
For the generating function $R(x)=\left[\frac{1}{1-t}\right]^{\alpha}$ the coefficients are determined by the expression:

$$
{n+\alpha-1 \choose n}.
$$
Then, according to (\ref{P_n_common}), we obtain

$$C_n^{(\alpha)}(x)={{\sum_{k=\left \lceil {{n}\over{2}} \right \rceil}^{n}{{{k
 }\choose{n-k}}\,{{k+\alpha-1\choose{k}}\,\left(-1\right)^{n-k}\,(2x)^{2\,k-n
 }}}}}.
$$ 

\section{Hermite Polynomials}

Hermite polynomials are a classical orthogonal polynomials that arise in probability, in combinatorics,  in numerical analysis, in finite element methods and in physics. They are also used in systems theory in connection with nonlinear operations on Gaussian noise.

The generating function for the Hermite polynomials  is given as follows \cite{Roman1984}:
$$
e^{2xt-t^2}=\sum_{n\geq 0} H_n(x)\frac{t^n}{n!}.
$$ 
Let us consider the generating function as a composition of generating functions
$
e^{f(x,t)}.
$ 

Then, according to (\ref{P_n_common}), we obtain
$$
H_n(x)=n!\sum_{k=\lceil {{n}\over{2}} \rceil}^n \frac{2^{2k-n}(-1)^{n-k}x^{2k-n}}{(n-k)!(2k-n)!}.
$$

\section{Associated Laguerre Polynomials}
Associated Laguerre polynomials are solutions of the equation:
$$
xy''+(\alpha+1-x)y'+ny=0.
$$

Laguerre polynomials arise in quantum mechanics, in the radial part of the solution of the Schroedinger equation for a one-electron atom.

The generating function for the Laguerre polynomials  is given as follows \cite{MathWolfram}:
$$
(1-t)^{-\alpha-1} e^{\frac{xt}{t-1}}=\sum_{n\geq 0} L_n^{(\alpha)}(x).
$$

Let us obtain the expression of the coefficients  for the generating function $(1-t)^{-\alpha-1}$

$$
L_1(n,\alpha)=\frac{1}{n!}{{\sum_{k=0}^{n}{\left(\alpha+1\right)^{k}\,\left(-1\right)^{n-k}\, \left[{n \atop k}\right]}}}
$$
or
$$
L_1(n,\alpha)={n+\alpha\choose n}.
$$

The expression of the coefficients for the generating function $e^{\frac{xt}{t-1}}$ is 
$$L_2(n,x)=\left\{
\begin{array}{ll}
\sum_{k=1}^{n}{{{\frac{(-1)^{k}}{k!}{n-1\choose k-1}\,x^{k
 }}}}, & n>0,\\
 1, & n=0.
\end{array} 
\right.$$ 

Therefore,

$$
L_n^{(\alpha)}(x)=\sum_{i=0}^n L_1(i,\alpha)L_2(n-i,x).
$$

\section{Stirling Polynomials}
The generating function for the Stirling polynomials  is given as follows \cite{MathWolfram}:

$$
\sum_{n\geq 0} S_n(x)\frac{t^n}{n!}=\left(\frac{t}{1-e^{-t}}\right)^{x+1}.
$$

Let us consider the generating function as a composition of generating functions

$$
\left(\frac{t}{1-e^{-t}}\right)^{x+1}=e^{-(x+1)\ln\left(1+\frac{e^{-t}-1}{-t}-1\right)}.
$$

The composita of the generating function $\frac{e^{-t}-1}{-t}-1$ is given as follows:

$$(-1)^n\sum_{j=0}^{k} \frac{j!\,\left(-1\right)^{k+j}\,{{k}\choose{j}}\,\left\{{n+j \atop j}\right\}}{\left(n+j\right)!}.$$

According to \cite{KruCompositae}, the composita of the generating function $-(x+1)\ln\left(1+\frac{e^{-t}-1}{-t}-1\right)$ is given as follows:

$$(-1)^{n+m}\,\left(\sum_{k=m}^{n}\frac{m!}{k!}{{{\left[{k \atop m}\right]\,\sum_{j=0}^{k}{{{\frac{j!(-1)^{k+j}}{(n+j)!}\,{{k}\choose{j}}\,\left\{{n+j \atop j}\right\}}}}}}}\right)\,\left(x+1\right)^{m}.$$

Therefore, we obtain

$$S_n(x)=\sum_{m=0}^{n}(-1)^{n+m}\,\left(\sum_{k=m}^{n}\frac{n!}{k!}{{{\left[{k \atop m}\right]\,\sum_{j=0}^{k}{{{\frac{j!(-1)^{k+j}}{(n+j)!}\,{{k}\choose{j}}\,\left\{{n+j \atop j}\right\}}}}}}}\right)\,\left(x+1\right)^{m}.
 $$

Let us consider other composition of generating functions
$$
\left(\frac{t}{1-e^{-t}}\right)^{x+1}=\left[\frac{1}{1+\frac{e^{-t}-1}{-t}-1}\right]^{x+1}.
$$

For the generating function $R(x)=\left[\frac{1}{1+t}\right]^{x+1}$ the coefficients are determined by the expression \cite{KruCompositae}:

$$
{n+x \choose n}(-1)^n.
$$
Then, according to (\ref{composition}), we obtain

$$S_n(x)=n!\sum_{k=0}^n{k+x \choose k}\,\sum_{j=0}^{k}{{{\frac{j!(-1)^{n+j}}{(n+j)!}\,{{k}\choose{j}}\, \left\{{n+j \atop j}\right\}}}}.
$$

\section{Abel Polynomials}
The Abel polynomials form a polynomial sequence of binomial type.

The generating function for the Abel polynomials  is given as follows \cite{MathWolfram}:

$$
e^{\frac{x}{a}W(at)}=\sum_{n\geqslant 0} A_n(x,a)\frac{t^n}{n!},
$$
where $W(t)$ is the Lambert W-function. The composita of $W(t)$ is given as follows
$$
W^{\Delta}(n,k)={{k\,n^{n-k-1}(-1)^{n-k}}\over{\left(n-k\right)!}}
$$

Then, according to (\ref{composition}), we obtain
$$A_n(x,a)=\sum_{k=1}^{n}{a^{n-k}\,k\,n^{n-k-1}\,\left(-1\right)^{n-k}\,
 {{n}\choose{k}}}\,x^{k}=x\sum_{k=0}^{n-1}{n-1 \choose k}x^k(-na)^{n-k-1}=x(x-na)^{n-1}.$$

\section{Bernoulli Polynomials of the Second Kind}
The Bernoulli polynomials are important in various branches of mathematics such as number theory, analysis, calculus of finite differences, statistics, etc.

The generating function for the Bernoulli Polynomials of the Second Kind is given as follows \cite{MathWolfram}:
$$
\frac{t(t+1)^x}{\ln(t+1)}=\sum_{n\geqslant 0} \frac{b_n(x)}{n!}t^n.
$$

First, we shall find an expression of the coefficients for the generating function  $\frac{t}{\ln(1+t)}=\frac{e^{\ln(t+1)}-1}{\ln(1+t)}$.
The composition $g(\ln(1+t))$ such that $g(t)=\frac{e^t-1}{t}$ defines the required expression of coefficients:
$$
\frac{1}{n!}\sum_{k=0}^n\frac{1}{k+1}\genfrac{[}{]}{0pt}{}{n}{k}
$$

Then,  we obtain
 
$$b_n(x)=\sum_{i=0}^{n}{n\choose i}\sum_{k=0}^{n-i}\frac{1}{k+1}\genfrac{[}{]}{0pt}{}{n-i}{k}\,\sum_{k=0}^{i}\genfrac{[}{]}{0pt}{}{i}{k}\,x^{k}$$  

or

$$b_n(x)=\sum_{i=0}^{n}{\frac{1}{(n-i)!}}\sum_{k=0}^{n-i}\frac{1}{k+1}{\genfrac{[}{]}{0pt}{}{n-i}{k}\,{x \choose i}}$$  

Let us consider other form of the generating function $$\frac{t}{\ln(1+t)}=\frac{1}{1+\frac{\ln(1+t)}{t}-1}.$$

The composita of the generating function $\frac{\ln(1+t)}{t}-1$ is given as follows:
$$
\sum_{j=0}^k\frac{j!}{(n+j)!} \genfrac{[}{]}{0pt}{}{n+j}{j}{k \choose j}(-1)^{k-j}
$$
Then, the coefficients of the generating function $\frac{t}{\ln(1+t)}$ are determined by the expression:
$$\sum_{k=0}^{n}{\sum_{j=0}^{k}{{{\left(-1\right)^{j}\,j!\,{{k
 }\choose{j}}\,\genfrac{[}{]}{0pt}{}{n+j}{j}}\over{\left(n+j
 \right)!}}}}$$

Therefore,  we obtain

$$b_n(x)=\sum_{i=0}^{n}\frac{1}{(n-i)!}{{{\left(\sum_{k=0}^{i}{\sum_{j=0}^{k}{{{\left(-1
 \right)^{j}\,j!\,\genfrac{[}{]}{0pt}{}{j+i}{j}\,{{k}\choose{j
 }}}\over{\left(j+i\right)!}}}}\right)\,\sum_{k=0}^{n-i}{
 \genfrac{[}{]}{0pt}{}{n-j}{k}\,x^{k}}}}}$$

\section{Generalized Bernoulli Polynomials}
Bernoulli polynomials play an important role in various expansions and
approximation formulas which are useful both in analytic theory of numbers
and in classical and numerical analysis.

The generating function for the Generalized Bernoulli polynomials is given as follows \cite{MathWolfram}:
$$
e^{xt}\left( \frac{t}{e^t-1}\right)^{\alpha}=\sum_{n\geqslant 0} B_n^{(\alpha)}(x)\frac{t^n}{n!}
$$

First, we shall get an expression of the coefficients for the generating function 
$$
\left( \frac{t}{e^t-1}\right)^{\alpha}=e^{-\alpha\ln\left(1+\frac{e^t-1}{t}-1\right)}
$$ 

The composita of the generating function  $\left(\frac{e^t-1}{t}-1\right)$ equals to

$$\sum_{j=0}^{k}{{{j!\,\left(-1\right)^{k+j}\,{{k}\choose{j}}\,
\genfrac{\{}{\}}{0pt}{}{n+j}{j}}\over{\left(n+j\right)!}}}$$
 Then the composita of the generating function  $-\alpha\ln\left(1+\frac{e^t-1}{t}-1\right)$ is
$$\alpha^{m}\,\left(-1\right)^{m}\,\sum_{k=m}^{n}\frac{m!}{k!}{{{
\genfrac{[}{]}{0pt}{}{k}{m}\,\sum_{j=0}^{k}{{{\frac{j!\,(-1
 )^{k+j}}{(n+j)!}\,{{k}\choose{j}}\,\genfrac{\{}{\}}{0pt}{}{n+j}{j}
 }}}}}}$$

Then, the coefficients of the generating function $\left( \frac{t}{e^t-1}\right)^{\alpha}$  are determined by the expression:
$$T_n(\alpha)=\sum_{m=0}^n\alpha^{m}\,\left(-1\right)^{m}\,\sum_{k=m}^{n}\frac{1}{k!}{{{
 \genfrac{[}{]}{0pt}{}{k}{m}\,\sum_{j=0}^{k}{{{\frac{j!\,(-1
 )^{k+j}}{(n+j)!}\,{{k}\choose{j}}\,\genfrac{\{}{\}}{0pt}{}{n+j}{j}
 }}}}}}$$

On the other hand we can write
$$
\left( \frac{t}{e^t-1}\right)^{\alpha}=\left[\frac{1}{1+\frac{e^{t}-1}{t}-1}\right]^{\alpha}=g(h(t)),
$$
where $g(t)=\left[\frac{1}{1+t}\right]^{\alpha}$ and $h(x)=\frac{e^t-1}{t}-1$.

Then, according to (\ref{composition}), we obtain
$$
T_n(\alpha)=\sum_{k=0}^n (-1)^k{k+\alpha-1 \choose \alpha-1}\sum_{j=0}^{k}{{{j!\,\left(-1\right)^{k+j}\,{{k}\choose{j}}\,
 \genfrac{\{}{\}}{0pt}{}{n+j}{j}}\over{\left(n+j\right)!}}}.
$$

Therefore, we get the explicit formula for the Generalized Bernoulli polynomials
$$
B_n^{(\alpha)}(x)=n!\sum_{i=0}^{n} T_i(\alpha)\frac{x^{n-i}}{(n-i)!}.
$$

\section{Euler Polynomials}

The generating function for the Euler Polynomials is given as follows \cite{MathWolfram}:
$$
\frac{2e^{xt}}{1+e^t}=\sum_{n\geqslant 0} E_n(x)\frac{t^n}{n!}.
$$
Represent this generating function to the following form:
$$
\frac{2e^{xt}}{1+e^t}=e^{xt}\frac{1}{1+\frac{1}{2}(e^t-1)}.
$$

Let us consider the following generating function  as a composition of generating functions 
$$f(g(t))=\frac{1}{1+\frac{1}{2}(e^t-1)},$$
where $f(t)=\frac{1}{1+t}$ and $g(t)=\frac{1}{2}(e^t-1)$.

Then, according to (\ref{composition}), the coefficients of the generating function $f(g(t))=\frac{1}{1+\frac{1}{2}(e^t-1)}$  are determined by the expression:

$${{\frac{1}{n!}\sum_{k=0}^{n}{{{{\frac{(-1)^{k}k!}{2^{k}}\genfrac{\{}{\}}{0pt}{}{n}{k}}}}}}}.
  $$ 

Therefore,according to the rule of multiplication of exponential generating functions, we obtain
$$
E_n(x)=\sum_{i=0}^n {n \choose i}x^i{{\sum_{k=0}^{n-i}{{{\frac{(-1)^{k}\,k!}{2^{k}}\,\genfrac{\{}{\}}{0pt}{}{n-i}{k}}}}}}.
$$  

\section{Peters Polynomials}

The Peters polynomials are a generalization of the Boole polynomials.

The generating function for the Peters polynomials is given as follows \cite{Roman1984}:

$$
\left[1+(1+t)^{\lambda}\right]^{-\mu}(1+t)^x=\sum_{n\geqslant 0} s_n(x,\lambda,\mu)\frac{t^n}{n!}.
$$

Represent the generating function $\left[1+(1+t)^{\lambda}\right]^{-\mu}$ as a composition of generating functions $e^t-1$ and $\ln(1+t)$:
$$
\left[1+(1+t)^{\lambda}\right]^{-\mu}=2^{-\mu}e^{-\mu\ln\left({1+\frac{e^{\lambda\ln(1+t)}-1}{2}}\right)}.
$$

The composita of the generating function
 $\left({\frac{e^{\lambda\ln(1+t)}-1}{2}}\right)$ is equal to
$$
{{\frac{m!}{2^{m}\,n!}\,\sum_{k=m}^{n}\genfrac{\{}{\}}{0pt}{}{k}{m}{\,
\genfrac{[}{]}{0pt}{}{n}{k}\,\lambda^{k}}}}.
$$

Then the composita of the generating function
${-\mu\ln\left({1+\frac{e^{\lambda\ln(1+t)}-1}{2}}\right)}$ is equal to the following expression:
$$
\frac{(-\mu)^ii!}{n!}\sum_{m=i}^n{{2^{-m}\,\sum_{k=m}^{n}{\genfrac{\{}{\}}{0pt}{}{k}{m}\,
 \genfrac{[}{]}{0pt}{}{n}{k}\,\lambda^{k}}}}\genfrac{[}{]}{0pt}{}{m}{i}.
$$

Hence, the coefficients of the generating function $\left[1+(1+t)^{\lambda}\right]^{-\mu}$ are determined by the expression:
$$
T_n(\lambda,\mu)=2^{-\mu}\sum_{i=1}^n\frac{(-\mu)^i}{n!}\sum_{m=i}^n{{2^{-m}\,\sum_{k=m}^{n}{\genfrac{\{}{\}}{0pt}{}{k}{m}\,
  \genfrac{[}{]}{0pt}{}{n}{k}\,\lambda^{k}}}}\genfrac{[}{]}{0pt}{}{m}{i}.
$$

Therefore, according to the rule of multiplication of exponential generating functions, we obtain
$$
s_n(x,\lambda,\mu)=2^{-\mu}\sum_{j=0}^n {n\choose j}\left(\sum_{i=1}^j(-\mu)^i\sum_{m=i}^n{{2^{-m}\,\sum_{k=m}^{j}{\genfrac{\{}{\}}{0pt}{}{k}{m}\, \genfrac{[}{]}{0pt}{}{j}{k}\,\lambda^{k}}}}\genfrac{[}{]}{0pt}{}{m}{i}\right)\sum_{k=0}^{n-j}  \genfrac{[}{]}{0pt}{}{n-j}{k}x^k 
$$

Let us consider other composition of generating functions
$$
[1+(1+t)^{\lambda}]^{-\mu}=\frac{1}{2^{\mu}}\left[\frac{1}{1+\frac{(1+t)^{\lambda}-1}{2}}\right]^{\mu}
$$

For the generating function $\frac{(1+t)^{\lambda}-1}{2}$ the composita is
$$
\frac{1}{2^k}\sum_{j=0}^k {k\choose j} {j\lambda\choose n}(-1)^{k-j}
$$

Then the coefficients of the composition of generating functions are determined by the expression:
$$T_n(\lambda,\mu)={{\frac{1}{2^{\mu}}\sum_{k=0}^{n}{{{\frac{1}{2^{k}}{{\mu+k-1}\choose{k}}\,\sum_{j=0}^{k}{\left(-1
 \right)^{j}\,{{k}\choose{j}}\,{{j\,\lambda}\choose{n}}}}.
 }}}}$$
 
Therefore, we obtain
$$
s_n(x,\lambda,\mu)=n!\sum_{i=0}^n { T_i(\lambda,\mu){x \choose n-i}}.
$$ 

The case for $\lambda=2$ and $\mu=-2$  is considered as a triangle of coefficients of degrees $x$ in the sequence $A137393$ \cite{oeis}.

\section{Narumi Polynomials}

The generating function for the Narumi Polynomials  is given as follows \cite{MathWolfram}:
$$
\left[\frac{t}{\ln(1+t)}\right]^{\alpha}(1+t)^x=\sum_{n\geq 0} S_n(x,\alpha)\frac{t^n}{n!}.
$$

Let us consider the left part of the generating function as a composition of generating functions $e^x$ and $\ln(1+x)$
$$
\left[\frac{t}{\ln(1+t)}\right]^{\alpha}=e^{-\alpha\ln\left(1+\frac{\ln(1+t)}{t}-1\right)}.
$$

The composita of the generating function $\left(\frac{\ln(1+t)}{t}-1\right)$ is given as follows:
$$
\left(\frac{\ln(1+t)}{t}-1\right)^k=\sum_{j=0}^k {k\choose j}\left(\frac{\ln(1+t)}{t}\right)^j(-1)^{k-j}.
$$

For the generating function $\left(\frac{\ln(1+t)}{t}\right)^j$ the coefficients are determined by the expression \cite{KruCompositae}:

$$
\left[{n+j \atop j}\right]\frac{j!}{(n+j)!}.
$$
Therefore, we obtain
$$
\left(\frac{\ln(1+t)}{t}-1\right)^k=\sum_{j=0}^k {k\choose j}\left[{n+j \atop j}\right]\frac{j!}{(n+j)!}(-1)^{k-j}t^j.
$$

Then the composita of the generating function  $-\alpha\ln\left(1+\frac{\ln(1+t)}{t}-1\right)$ is 
$$
(-\alpha)^m\sum_{k=m}^n\sum_{j=0}^k {k\choose j}\left[{n+j \atop j}\right]\frac{j!}{(n+j)!}(-1)^{k-j}\left[{k \atop m}\right]\frac{m!}{k!}.
$$

For the generating function $e^{-\alpha\ln\left(1+\frac{\ln(1+t)}{t}-1\right)}$ the coefficients are determined by the expression 
$$
\sum_{m=0}^n(-\alpha)^m\sum_{k=m}^n\sum_{j=0}^k {k\choose j}\left[{n+j \atop j}\right]\frac{j!}{(n+j)!}(-1)^{k-j}\left[{k \atop m}\right]\frac{m!}{k!}\frac{1}{m!}.
$$
After transformation we obtain
$$
\sum_{m=0}^n(-\alpha)^m\sum_{k=m}^n\sum_{j=0}^k \left[{n+j \atop j}\right]\frac{(-1)^{k-j}}{(k-j)!(n+j)!}\left[{k \atop m}\right].
$$

Therefore, according to the rule of a product of exponential generating functions, we obtain 
$$S_n(x,\alpha)=\sum_{i=0}^{n}{\left(i!\sum_{m=0}^{i}{(-\alpha)^{m}\,\left(\sum_{k=m
 }^{i}{\left[{k \atop m}\right]\,\sum_{j=0}^{k}{{{
 \left[{i+j \atop i}\right]\,\left(-1\right)^{j-k}}\over{
 \left(j+i\right)!\,\left(k-j\right)!}}}}\right)
 }\right)\,{{n}\choose{i}}\,\sum_{k=0}^{n-i}{\left[{n-i \atop k}\right]\,x^{k}}}.$$

Let us consider other composition of generating functions
$$
\left[\frac{1}{1+\frac{\ln(1+t)}{t}-1}\right]^{\alpha}.
$$
The coefficients are determined by the expression
$$
\sum_{k=0}^n{k+\alpha-1 \choose k}(-1)^k\sum_{j=0}^k {k\choose j}\left[{n+j \atop j}\right]\frac{j!}{(n+j)!}(-1)^{k-j}.
$$
Therefore,
$$S_n(x,\alpha)=n!\sum_{i=0}^{n}{\left(\sum_{k=0}^{i}{{{k+\alpha-1}\choose{k}}
 \,\sum_{j=0}^{k}{{{\left(-1\right)^{j}\,j!\,\left[{i+j \atop j}\right]\,{{k}\choose{j}}}\over{\left(j+i\right)!}}}}\right)\,{{
 x}\choose{n-i}}}.$$
 
\section{Humbert Polynomials}

The generating function for the Humbert polynomials  is given as follows \cite{PolyExpan}:
$$
\left(1-3xt+t^3\right)^{-\lambda}
$$

The composita of the generating function $-3xt+t^3$  is given as follows
$$\frac{1}{2}{{{{k}\choose{{{n-k}\over{2}}}}\,3^{{{3\,k-n}\over{2}}}\,\left(
 \left(-1\right)^{n-k}+1\right)\,\left(-x\right)^{{{3\,k-n}\over{2}}}
 }}.$$

Then, the coefficients of the generating function $\exp\left(-\lambda\ln\left(1-3xt+t^3\right)\right)$ are determined by the expression:
$$H_n(x,\lambda)=\sum_{m=0}^{n}{\left(-1\right)^{m}\,\left(\sum_{j=0}^{{{n-m}\over{2
  }}}{{{\frac{3^{n-3\,j}}{(n-2j)!}\,{{n-2\,j}\choose{j}}\,\genfrac{[}{]}{0pt}{}{n-2j}{m}\,\left(-x\right)^{n-3\,j}}}}
 \right)\,\lambda^{m}}$$

Let us consider the generating function as the following composition of generating functions 
$$
\left[\frac{1}{1-h(t)}\right]^{\lambda},
$$
where $h(t)=3xt-t^3$.

Then
$$
H_n(x,\lambda)=(-1)^n\sum_{k=1}^n \frac{1}{2}{{{{k}\choose{{{n-k}\over{2}}}}\,3^{{{3\,k-n}\over{2}}}\,\left(
 \left(-1\right)^{n-k}+1\right)\,\left(-x\right)^{{{3\,k-n}\over{2}}}
 }} {k+\lambda-1 \choose k}.
$$

Substituting $2m$ for $n-k$, we get the following expression
$$
H_n(x,\lambda)=\sum_{m=0}^{\frac{n}{3}} {{{{n-2m}\choose{{m}}}\,3^{{{n-3m}}}\,\,(-1)^m\,x^{{{n-3m}}} }} {n-2m+\lambda-1 \choose n-2m}.
$$

\section{Lerch Polynomials}

The generating function for the Lerch Polynomials  is given as follows \cite{PolyExpan}:
$$
(1-x\ln(1+t))^{\lambda}=\sum_{n\geqslant 0} L_n(x,\lambda)
$$

Since the composita of $\ln(1+t)$ equals to $\genfrac{[}{]}{0pt}{}{n}{k}\frac{k!}{n!}$ and the coefficients of the generating function $\left[\frac{1}{1-t}\right]^{\lambda}$ are determined by ${n+\lambda-1 \choose \lambda-1}$, we obtain

$$L_n(x,\lambda)=\frac{1}{n!}{{\sum_{k=1}^{n}{k!\,\genfrac{[}{]}{0pt}{}{n}{k}\,x^{k}\,{{
 \lambda+k-1}\choose{k}}}}}$$

\section{Mahler Polynomials}

The generating function for the Mahler Polynomials  is given as follows \cite{Roman1984}:
$$
e^{x(1+t-e^t)}=\sum_{n\geqslant 0} M_n(x)\frac{t^n}{n!}
$$

The composita of the sum of the generating functions $G(t)=t$ and $R(t)=-(e^t-1)$ is equal to
$$\sum_{j=0}^{k}{{{\frac{j!}{(n-k+j)!}\,\left(-1\right)^{k-j}\,{{k}\choose{j}}\,
\genfrac{\{}{\}}{0pt}{}{n-k+j}{j} }}}.
 $$

Therefore, according to (\ref{composition}, we obtain
$$M_n(x)=\sum_{k=0}^{n} x^{k}{\left(\sum_{j=0}^{k}{\left(-1\right)^{j}\,{{n
 }\choose{k-j}}\,\genfrac{\{}{\}}{0pt}{}{n-k+j}{j}}\right)}$$ 

For instance, the coefficients of Mahler polynomials
is considered as a triangle of coefficients of degrees $x$ in the sequence $A137375$ \cite{oeis}.

\end{document}